\documentclass[reqno, 11pt]{amsart}

\newtheorem{theorem}{Theorem}[section]
\newtheorem{lemma}[theorem]{Lemma}
\newtheorem{proposition}[theorem]{Proposition}
\newtheorem{corollary}[theorem]{Corollary}

\theoremstyle{definition}
\newtheorem{definition}[theorem]{Definition}
\newtheorem{ex}[theorem]{Example}

\theoremstyle{remark}
\newtheorem{remark}[theorem]{Remark}

\numberwithin{equation}{section}

\usepackage{bbm}
\usepackage{url}
\usepackage{euscript}
\usepackage{pb-diagram}
\usepackage{lamsarrow}
\usepackage{pb-lams}
\usepackage{amsmath}
\usepackage{amsthm}

\usepackage{epsfig}
\usepackage{amssymb}
\usepackage{pifont}
\usepackage{amsbsy}

\newskip\aline \newskip\halfaline
\def\skipaline{\vskip\aline}

\def\qedbox{$\rlap{$\sqcap$}\sqcup$}
\def\qed{\nobreak\hfill\penalty250 \hbox{}\nobreak\hfill\qedbox\skipaline}

\def\proofend{\eqno{\mbox{\qedbox}}}

\newcommand{\one}{\mathbbm{1}}

\newcommand\bC{{\mathbb C}}

\newcommand\bR{{\mathbb R}}

\newcommand\bZ{{\mathbb Z}}

\newcommand{\bm}{\boldsymbol{m}}

\newcommand{\bsM}{\boldsymbol{M}}
\newcommand{\bsN}{\boldsymbol{N}}

\newcommand{\bsS}{\boldsymbol{S}}

\newcommand{\bdel}{\boldsymbol{\delta}}
\newcommand{\beps}{\boldsymbol{\epsilon}}
\newcommand{\blam}{\boldsymbol{\lambda}}
\newcommand{\bmu}{\boldsymbol{\mu}}
\newcommand{\bpi}{\boldsymbol{\pi}}
\newcommand{\brho}{\boldsymbol{\rho}}

\newcommand{\bsi}{\boldsymbol{\sigma}}

\newcommand{\eps}{{\epsilon}}
\newcommand{\vfi}{{\varphi}}

\newcommand{\eB}{\EuScript{B}}

\newcommand{\eD}{\EuScript{D}}

\newcommand{\eI}{\EuScript{I}}

\newcommand{\eL}{\EuScript{L}}

\newcommand{\eN}{\EuScript{N}}

\newcommand{\eS}{\EuScript{S}}

\newcommand{\eW}{\EuScript{W}}

\newcommand{\ra}{\rightarrow}
\newcommand{\hra}{\hookrightarrow}
\newcommand{\Lra}{{\longrightarrow}}

\def\inpr{\mathbin{\hbox to 6pt{\vrule height0.4pt width5pt depth0pt \kern-.4pt \vrule height6pt width0.4pt depth0pt\hss}}}

\DeclareMathOperator{\Cr}{\mathbf{Cr}}

\DeclareMathOperator{\Lag}{Lag_h}

\newcommand{\whbm}{\widehat{\boldsymbol{m}}}

\begin{document}

\title{On a Bruhat-like poset}
\date{Started October 24, 2007. Completed November 2, 2007.}

\author{Liviu I. Nicolaescu}

\address{Department of Mathematics, University of Notre Dame, Notre Dame, IN 46556-4618.}
\email{nicolaescu.1@nd.edu}

\begin{abstract}  We  investigate the poset of strata of a Schubert like stratification  on the Grassmannian of   hermitian lagrangian spaces in $\bC^n\oplus \bC^n$. We prove that this poset is a modular lattice, we compute  its M\"{o}bius function and we investigate its order intervals.

\end{abstract}

\maketitle

\tableofcontents

\section*{Introduction}
\addtocontents{section}{Introduction}
\setcounter{equation}{0}
The   Grassmannian of the  hermitian lagrangian subspaces of
$E=\bC^n\oplus\bC^n$ is the space $\Lag(n)$ of  complex subspaces $L \subset E$ such that
\[
L^\perp =JL\;\;\mbox{where}\;\;L= \left[\begin{array}{lc}
0 &-\one_{\bC^n}\\
\one_{\bC^n} & 0
\end{array}
\right].
\]
This space is naturally  a submanifold of the space  $\eS_E$ of   selfadjoint operators $E\ra E$. As such, it is equipped with a natural Riemann metric. In \cite{N} we proved that a certain  linear   function $ f:\eS_E\ra \bR$  induces   a perfect Morse function on $\Lag(n)$. Moreover, there exists a bijection  between the collection $2^{[n]}$ of subsets of $[n]=\{1,\dotsc, n\}$ and $\Cr_f$, the set of critical points of $f$,
\[
2^{[n]}\ni I\mapsto L_I\in \Cr_f.
\]
The  negative gradient  flow of $f$ satisfies the Smale transversality condition. If $W_I^-$ denotes the unstable manifold of $L_I$ then the collection
\[
\eW_n:=\bigl\{ W_I^-;\;\;I\in 2^{[n]}\,\bigr\}
\]
defines a  Whitney regular stratification of $\Lag(n)$ which can be
given a Schubert-like description in terms of incidence conditions.

The collection $\eW_n$ of strata  is equipped with a Bruhat-like partial order $\prec$  defined by
\[
W_I^-\prec W_J^- \Longleftrightarrow  W_I^-\subset {\rm closure}\,(W_I^-).
\]
More explicitly,
\[
W_I^-\prec W_J^- \Longleftrightarrow \# (I\cap [k,n]) \leq \# (J\cap [k,n]),\;\;\forall k=1,\dotsc, n.
\]
 We obtain in this fashion a partial order $<_n$ on $2^{[n]}$ and we denote the resulting poset by $\eL_n$. The Hasse diagrams of $\eL_n$, $1\leq n\leq 4$ are depicted in Figure \ref{fig: mob1}.

It is perhaps  instructive to give another description of this partial order in terms of the game  \emph{beads-along-a-rod}.

Suppose   are given a thin  rod with   points marked $1,\dotsc, n$ in increasing linear order, left-to-right along the rod.     To any  subset $I\subset [n]$ we associate a configuration of beads along this rod placed in the positions $i\in I$.  The game consist of a succession of elementary moves of two types.

\begin{itemize}
\item  Slide one bead to the  left neighboring position, if that position is unoccupied by another ring.

\item   Remove  the  bead  on  the leftmost position $1$,  if  there is such a bead.

\end{itemize}
We  will refer to these moves as \emph{elementary left slides}.  If $I, J\in 2^{[n]}$ describe two configurations of rings, then $I\prec J$ if one can go from the configuration $J$  to the configuration $I$ by a sequence of elementary  left slides.

The poset $\eL_n$ is closely related to the Boolean poset  $\eB_n$ of subsets of $[n]$ ordered by inclusion.    To describe this relationship we need to go back to Morse theory.

We denote by $W_I^+$ the \emph{stable} manifold of the critical point $L_I\in \Lag(n)$. Then
\begin{equation*}
W_J^-\prec  W_I^- \Longleftrightarrow W_J^+ \cap W_I^-\neq\emptyset.
\tag{G}
\label{tag: g}
\end{equation*}
In \cite{N} we have proved that the stable/unstable  manifolds  $W^\pm_I$ determine homology classes $[W_I^\pm]\in H_\bullet(\Lag(n),\bZ)$. The    partial  order  on $\eB_n$ can also be formulated  in terms  of \emph{homological} intersection. More precisely
\begin{equation*}
J\subset I \Longleftrightarrow [W_J^+]\bullet [W_I^-]\neq 0,
\tag{H}
\label{tag: h}
\end{equation*}
where $\bullet$ denotes the intersection  pairing in homology.

On complex   Grassmannians   the two statements  (\ref{tag: g}) and (\ref{tag: h}) are equivalent and they   define the Bruhat order on the set of Schubert  varieties. On the Grassmannian $\Lag(n)$ we  only have the implication (\ref{tag: h}) $\Rightarrow$ (\ref{tag: g}). In particular,  the  tautological map $\eB_n\ra \eL_n$ is increasing, but its inverse is not.

The poset  $\eL_n$  shares many combinatorila features with the Bruhat posets. We show (Proposition \ref{prop: lattice}, \ref{prop: modular}) that $\eL_n$ is a \emph{modular ortholattice  lattice} with rank function $\brho:\eL_n\ra \bZ$ given by
\[
\brho(I)=\sum_{i\in I} i,
\]
and the complement map $\bsi_n:\eL_n\ra \eL_n$ given by $S\mapsto \{1,\dotsc,n\}\setminus S$.  As explained   in \cite{Wach}  the modularity implies that the poset $\eL_n$ is homotopy Cohen-Macaulay.

We define a pair $(I,J)\subset \eL_n\times \eL_n$ to be \emph{elementary} if the sequence of  integers
\[
\delta_k= \# \bigl(\,J\cap [k,n]\,\bigr)- \# \bigl(\,I\cap [k,n]\,\bigr), \;\;k=1,2,\dotsc, n
\]
consists only of $0$'s and $1$'s   and   there are no consecutive $1$'s.  If $\bmu_n$ denotes the M\"{o}bius function of $\eL_n$, then   we show (Theorem \ref{th: mob}) that
\[
\bmu_n(S,T)=\begin{cases}
(-1)^{\brho(J)-\brho(I)} & \mbox{$(I,J)$ is  an elementary pair},\\
0 & \mbox{otherwise}.
\end{cases}
\]
Since a modular lattice  is shellable we deduce from \cite[Thm. 5.6]{BjWa}  that if $(I,J)$ is not an elementary pair then the   open order interval $(I,J)_{\eL_n}$ is a contractible poset.  In Proposition \ref{prop: ball}  we prove a stronger state namely that in this case the   nerve of the open order interval $(I,J)_{\eL_n}$  is homeomorphic to the closed Euclidean ball of dimension $\brho(J)-\brho(I)-2$.

If $(I,J)$ is an elementary pair, then the shellability of $\eL_n$ together with the equality $\bmu_n(I,I)=(-1)^{\brho(J)-\brho(I)}$ imply via \cite[Thm. 5.6]{BjWa} that the open interval  $(I,J)_{\eL_n}$ is homotopic to  a sphere  of dimension $\brho(J)-\brho(I)-2$.  We were able to prove a slightly stronger result.  Namely,  we show that  if $(I,J)$ is an elementary pair then the order interval $[I,J]_{\eL_n}$ is isomorphic to the   boolean poset $\eB_{\brho(J)-\brho(I)}$ (Theorem \ref{th: boole}).

The key technical  device  that allowed us to reach these conclusions is  a  certain surgery-like operation $\#$ on pairs of  posets which we introduce  in  Section \ref{s: 1}. In this section we also describe explicitly the effect of    this surgery on the M\"{o}bius functions (Theorem \ref{th: mm}).

In  Section 2 we investigate a special case of this surgery operation  on a special class of posets we called \emph{layered}.     Given a layered poset $P$  we define its \emph{double}  as the poset  obtained by applying the  surgery   operation  in Section 1 to  two copies of $P$.  In  this section we analyze a  few special features of this doubling   operation.

In Section 3 we apply  the  general results to the poset $\eL_n$ which has  a natural    layer structure.   The key fact which allowed us to apply the general theory  developed in the previous  section is the ``surgery formula'' in  Proposition \ref{prop: double} which  states that $\eL_{n+1}$ is the double of the layered poset $\eL_n$.

\section{The M\"{o}bius functions of ``connected sums'' of  posets}
\label{s: 1}
\setcounter{equation}{0}
We follow  closely the poset terminology in \cite{Aig, Sta}. For simplicity we will concentrate exclusively on finite posets. If $(P, <)$ is such a poset,   we denote by $\bC^P$ the space of functions $P\ra \bC$, and we  define  a $\bZ$-linear map (integration)
\[
\bsS_P: \bC^P\ra \bC^P,\;\;  \bC^P\ni f\mapsto \eI_P f\in \bC^P,\;\;\eI_Pf(x)=\sum_{y\geq x} f(y),\forall x\in P.
\]
The  vector space $\bC^P$ has a canonical basis consisting of the Dirac functions
\[
\delta_x\in \bC^P,\;\;x\in P,\;\;\delta_x(y)=\begin{cases}
1 & y=x\\
0 & y\neq x.
\end{cases}
\]
With respect to this basis the above operator is described  by  its \emph{incidence matrix} $\eS_P\in \bC^{P\times P}$ defined by
\[
\eS_P(x,y)=\begin{cases}
1 & x\leq y\\
0 & x\nleqslant y.
\end{cases}
\]
The matrix $\eS_P$ is ``upper triangular'' and it can be   written as a sum $\eS_P=\one +\eN_P$, where  $\eN_P$ is ``strictly upper triangular'', i.e.,
\[
\eN_P(x,y)\neq 0\Longrightarrow x<y.
\]
This shows that the linear operator $\bsN_P$ determined by $\eN_P$ is nilpotent.  Hence $\bsS_P$ is invertible and its inverse is given by
\[
\bsM_P=\bsS_P^{-1}=\sum_{k\geq 0}(-1)^k \bsN_P^k.
\]
The matrix describing $\bsM_P$ in the Dirac basis is called the \emph{M\"{o}bius function of $P$} and it is denoted by $\bmu_P$. The equalities $\bsM_P\bsS_P=\one=\bsS_P^\dag\bsM_P^\dag$ (${}^\dag=$ transpose) translate into the recursion
\begin{equation}
\sum_{x\leq z\leq y}\bmu(x,z)=\delta_{x}(y)=\sum_{x\leq z\leq y}\bmu(z,y).
\label{eq: recurs}
\end{equation}
The M\"{o}bius function appears in the M\"{o}bius inversion formula
\begin{equation}
s(x)=\sum_{y\geq x}  f(y),\;\;\forall x\in P \Longleftrightarrow f(x)=\sum_{y\geq x} \bmu(x,y) s(y),\;\;\forall x\in P.
\label{eq: mob-for}
\end{equation}

Suppose we are given three posets $(Q, <_Q)$, $(P_k, <_k)$, $k=0,1$, and  injections $i_k: Q\ra P_k$ such that
\[
q <_Q q' \Longleftrightarrow i_k(q) <_k i_k(q'),\;\;\forall k=0,1.
\]
In other words, $Q$ is an induced subposet of both $P_0$ and $P_1$.  To simplify the presentation we will write $x\leq_0 q$ instead of $x\leq_0i_0(q)$, and $q\leq_1 y$ instead of $i_1(q)\leq_1 y$.

We  can  now   define a partial order $\prec$ on the disjoint sum $P_0\sqcup P_1$   by  setting $x\prec y$ if and only if

\begin{itemize}

\item either  both $x,y$ belong to the same set $P_k$ and $x<_ky$,

\item or $x\in P_0$, $y\in P_1$ and there exists $q\in Q$ such that $x\leq_0 q\leq_1 y$.

\end{itemize}

\begin{figure}[h]
\centerline{\epsfig{figure=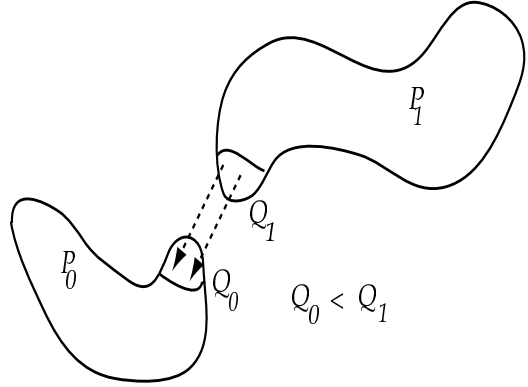, height=2in,width=2.7in}}
\caption{\sl Connecting $P_0$ to $P_1$ along $Q_0$ and $Q_1$.}
\label{fig: mob3}
\end{figure}

We denote this poset by $P_0 {}_{Q_0}\!\#_{Q_1}P_1$, where $Q_k=i_k(Q)\subset P_k$,  and we will refer to it as the \emph{connect sum of $P_0$ and $P_1$ along $Q_0$ and $Q_1$} (see Figure \ref{fig: mob3}).

\begin{ex}  Observe that if $P_0=P_1=Q$ and  $i_k=\one_Q$ then  the poset $P_0 {}_{Q_0}\!\#_{Q_1}P_1$ is isomorphic  with the poset  $\{-1,1\} \times P$  with the product  order
\[
(\eps_0,x_0)\leq (\eps_1,x_2)\Longleftrightarrow \eps_0\leq \eps_1,\;\;x_0\leq_Px_1.\proofend
\]
\end{ex}

\begin{proposition}  Let $P_0\stackrel{i_0}{\hookleftarrow}Q\stackrel{i_1}{\hra} P_1$  be  injective  increasing. We  set $Q_k:=i_k(Q)\subset P_k$, and form the connect sum $P=P_0 {}_{Q_0}\!\#_{Q_1}P_1$. Denote by $\eS_k$ and respectively $\bmu_k$ the incidence matrix and respectively the M\"{o}bius  function of $P_k$, $k=0,1$. If $x,y\in P_0\sqcup P_1$ then
\[
\bmu_{P}(x,y)=\begin{cases}
\bmu_k(x,y) & \mbox{if} \;\;  x,y\in P_k\\
-\sum_{x\leq_0 p_0\prec p_1\leq_1 y} \bmu_0(x,p_0)\bmu_1(p_1,y)& \mbox{if} \;\;x\in P_0,\;y\in P_1\\
0 & \mbox{otherwise}.
\end{cases}
\]
\label{prop: conn-sum}
\end{proposition}

\proof

With respect to the  direct sum decomposition $\bC^{P_0\sqcup P_0}=\bC^{P_0}\oplus \bC^{P_1}$ the  incidence matrix  $\eS_{P}$ has the block decomposition
\[
\eS_{P} =\left[\begin{array}{cc}
\eS_{0}  & \eB\\
0 &\eS_1
\end{array}
\right],\;\;\eS_k=\eS_{P_k},
\]
where $\eB: P_0\times P_1\ra \bZ$ satisfies
\[
\eB(x,y)=  \begin{cases}
 1 & \exists q\in Q:\;\; x\leq_0 q\leq_1 y\\
 0 & \mbox{otherwise}.
 \end{cases}
 \]
 Then
 \begin{equation}
 \bmu_{P}= \eS_{P}^{-1} = \left[\begin{array}{cc}
\eS_{0}^{-1}  & -\eS_0^{-1}\eB\eS_1^{-1}\\
0 &\eS_1^{-1}
\end{array}
\right]= \left[\begin{array}{cc}
\bmu_{0}  & -\bmu_0\eB\bmu_1\\
0 &\bmu_1
\end{array}
\right].
\label{eq: inverse}
\end{equation}
Note that if $x\in P_0$, $y\in P_1$ we then have
\[
\bmu_{P}(x,y)=-\sum_{p_0\in P_0,p_1\in P_1} \bmu_0(x,p_0)\eB(p_0,p_1)\bmu_1(p_1,y)=-\sum_{x\leq_0 p_0\prec p_1\leq_1 y} \bmu_0(x,p_0)\bmu_1(p_1,y).\proofend
\]

To  describe our next special case we need to introduce  some notation.  For any poset $P$ and any $X\subset P$ we set
\[
P^{\geq X}=\bigl\{ p\in P;\;\;\exists x\in X:\;\;p\geq x\,\bigr\},\;\;P^{\leq A}= \bigl\{ p\in P;\;\;\exists x\in X;\; p\leq x\,\bigr\}.
\]
Suppose now that  the induced subposet $Q_0\subset P_0$  satisfies the property
\begin{equation*}
\mbox{ $\forall x\in P_0$, $Q_0\cap P_0^{\geq x}$ is either empty or contains a unique minimal element $q_x$.}
\tag{$\bsM_+$}
\label{tag: m+}
\end{equation*}
and the subposet $Q_1\subset P_1$ satisfies
\begin{equation*}
\mbox{ $\forall y\in P_1$, $Q_1\cap P_1^{\leq y}$ is either empty or contains a unique maximal element $q^y$.}
\tag{$\bsM_-$}
\label{tag: m-}
\end{equation*}
Note that $\bsM_+$ implies that if $x\in P_0$ and $y\in Q_0$ then
\[
x\leq y \Longleftrightarrow q_x\leq y
\]
while $\bsM_-$ implies that if $x\in Q_1$ and $y\in P_1$ then
\[
x\leq y \Longleftrightarrow x\leq q^y.
\]

\begin{theorem} Suppose $Q_0\subset P_0$ satisfies (\ref{tag: m+}) and $Q_1\subset P_1$ satisfies (\ref{tag: m-}), and set $P= P_0\,{}_{Q_0}\!\#_{Q_1} P_1$,
\[
P_0^{\leq Q_0}=\bigl\{\, x\in P_0;\;\;P_0^{\geq x}\cap Q_0\neq \emptyset\,\bigr\},\;\; P_1^{\geq Q_1}= \bigl\{ \,y\in P_1;\;\;P_1^{\leq x}\cap Q_1\neq \emptyset\,\bigr\}.
\]
Then
\begin{equation}
\bmu_P(x,y) =\begin{cases}
\bmu_k(x,y) &\mbox{if}\;x,y\in P_k, \;k=0,1,\\
-\bmu_Q(q, q')  & \mbox{if} \;\;x=i_0(q)\in Q_0,\;y=i_1(q')\in Q_1,\\
0 &\mbox{otherwise}.
\end{cases}
\label{eq: mm}
\end{equation}
\label{th: mm}
\end{theorem}

\proof  The injections $i_0$ and $i_1$ define  injections $j_0=i_0\circ i_1^{-1}:Q_1\ra P_0$ and   $j_1:Q_0\ra P_1$. Let $x\in P_0^{\leq Q_0}$, $y\in P_1^{\geq Q_1}$. Then
\[
x\prec y \Longleftrightarrow  j_1(q_x) \leq q^y \leq y \Longleftrightarrow x\leq q_x\leq j_0(q^y),
\]
and
\begin{equation}
\eB(x,y)=\eS_1(j_1(q_x),y) = \eS_0(x, j_0(q^y)).
\label{eq: BS}
\end{equation}
In particular,
\begin{equation}
s\in P_0^{\leq Q_0},\;\; t \in P_1^{\geq Q_1},\;\;j_1(q_s)  \not\leq t \Longrightarrow \eB(s,t)=0.
\label{eq: vanish}
\end{equation}

Using (\ref{eq: inverse}) we deduce that if $x\prec y$ we have
\[
\bmu_P(x,y)=-\sum_{x\leq s\leq j_0(q^y),\;j_1(q_x)\leq t\leq y\;}\bmu_0(x,s)\eB(s,t)\bmu_1(t,y)
\]
\[
\stackrel{(\ref{eq: vanish})}{=}-\sum_{x\leq s\leq j_0(q^y)}\bmu_0(x,s)\Biggl( \sum_{j_1(q_s)\leq t\leq y\;} \eB(s,t)\bmu_1(t,y)\Biggr)
\]
\[
\stackrel{(\ref{eq: BS})}{=} -\sum_{x\leq s\leq j_0(q^y)}\bmu_0(x,s)\Biggl( \sum_{j_1(q_s)\leq t\leq y\;} \eS_1(j_1(q_s),t)\bmu_1(t,y)\Biggr)=  -\sum_{x\leq s\leq j_0(q^y)}\bmu_0(x,s)\delta_y(j_1(q_s)).
\]
The last sum is nonzero if and only if there exists $s\in P_0$ such that $x\leq s\leq j_0(q^y)$ and $y=j_1(q_s)\in Q_1$.  This can only happen when $y\in Q_1$, so that $y=q^y$, and $s$ is such that $j_1(q_s)=y$.

On the other hand, we have
\[
\bmu_P(x,y)\stackrel{(\ref{eq: vanish})}{=}-\sum_{j_1(q_x)\leq t\leq y}\Biggl(\sum_{x\leq s\leq j_0(q^t)}\bmu_0(x,s)\eB(s,t)\Biggr)\bmu_1(t,y)
\]
\[
\stackrel{(\ref{eq: BS})}{=} -\sum_{j_1(q_x)\leq t\leq y}\Biggl(\sum_{x\leq s\leq j_0(q^t)}\bmu_0(x,s)\eS_0(s,j_0(q^t))\Biggr)\bmu_1(t,y)=-\sum_{j_1(q_x)\leq t\leq y}\delta_{x}(j_0(q^t))\bmu_1(t,y).
\]
 We see that this sum can have a nontrivial term only if  $x\in Q_0$ so that $x=q_x$. Thus $x=q_0\in Q_0, \;\;y=q_1\in Q_1$. We have
\begin{equation}
\bmu_P(q_0,q_1)=-\sum_{j_1(q_0)\leq t\leq q_1}\delta_{j_1(q_0)}(q^t)\bmu_1(t, q_1)=-\sum_{q_0\leq s\leq j_0(q_1)}\bmu_0(q_0,s)\delta_{j_0(q_1)}(q_s)
\label{eq: mobq}
\end{equation}
For $q,q'\in Q$ we set
\[
A(q,q')= -\bmu_P(i_0(q),i_1(q')).
\]
From (\ref{eq: mobq}) we deduce
\[
A(q,q)=1,\;\; A(q,q')=\sum_{s\geq i_0(q), q_s=i_0(q')} \bmu_0(q,s).
\]
Hence
\[
\sum_{q\leq q'' \leq q'} A(q,q'')= \sum_{s\geq q,\; q_s\leq q''}\bmu_0(q,s)=\sum_{q\leq s\leq q'} \bmu_0(q,s)\stackrel{(\ref{eq: recurs})}{=}\delta_q(q'').
\]
We see that the function  $A: Q\times Q\ra \bZ$ satisfies the recurrence (\ref{eq: recurs}) so that
\[
A(q,q')= \bmu_Q(q,q').
\]
This completes the proof of Theorem \ref{th: mm}. \qed

\begin{remark} The conditions ($\bsM_\pm$) are not as restrictive as they look. For example,  they are automatically satisfied if the poset $Q$ is a lattice.\qed
\label{rem: mm}
\end{remark}

\begin{corollary}  If $P_0,P_1,Q$ are as in Theorem \ref{th: mm} then
\[
{\rm Range}\,\bmu_{P_0{}_{Q_0}\!\#_{Q_1} P_1}\subset \{0\}\cup {\rm Range}\,\bmu_{P_0}\cup {\rm Range}\,\bmu_{P_1}\cup {\rm Range}\,\bmu_{Q}.\proofend
\]
\label{cor: range}
\end{corollary}

If in Theorem \ref{th: mm} we choose $Q$ to  consists of a single point we obtain the following result.

\begin{corollary} Let $q_k\in  P_k$, $k=0,1$, and denote by $P$ the poset $P=P_0{}_{q_0}\!\#_{q_1}P_1$. Then
\begin{equation}
\bmu_P(x,y)= \begin{cases}
\bmu_k(x,y) & \mbox{if} \;\;  x,y\in P_k,\;k=0,1,\\
-\delta_{q_0}(x)\delta_{q_1}(y)& \mbox{if} \;\;x\in P_0,\;y\in P_1,\\
0 & \mbox{if} \;\; x\in P_1,\;\;y\in P_0.
\end{cases}
\label{eq: ex1}
\end{equation}
\end{corollary}

\section{Layered posets}
\label{s: 2}
\setcounter{equation}{0}
 We define a \emph{layer} structure  on a poset $P$ to be a pair  $(\beps,\blam)$ satisfying the following conditions.
 \begin{itemize}

 \item  $\beps$ is an increasing map $\beps: P\ra \{-1,1\}$ called the \emph{sign map} of the layer. We set $P^\pm:=\beps^{-1}(\pm)$. $P^-$ is called the \emph{lower layer} and $P^+$ is called the \emph{upper layer}.

 \item $\blam$ is a poset isomorphism $\blam:P^-\ra P^+$  such that $x<\blam(x)$, $\forall x\in P^-$. The map $\blam$ is called the \emph{lifting} map. Its inverse $\bdel:P^+\ra P^-$  is called the \emph{drop map}.

 \end{itemize}

 If $P$ is a lattice, we say that the layer structure is \emph{compatible with the lattice structure} if the layers are sublattices, i.e.,
 \[
 P^\pm\wedge P^\pm \subset P^\pm,\;\;P^\pm\vee P^\pm \subset P^\pm,
 \]
 and the  lifting map is an \emph{isomorphism} of lattices

\begin{figure}[h]
\centerline{\epsfig{figure=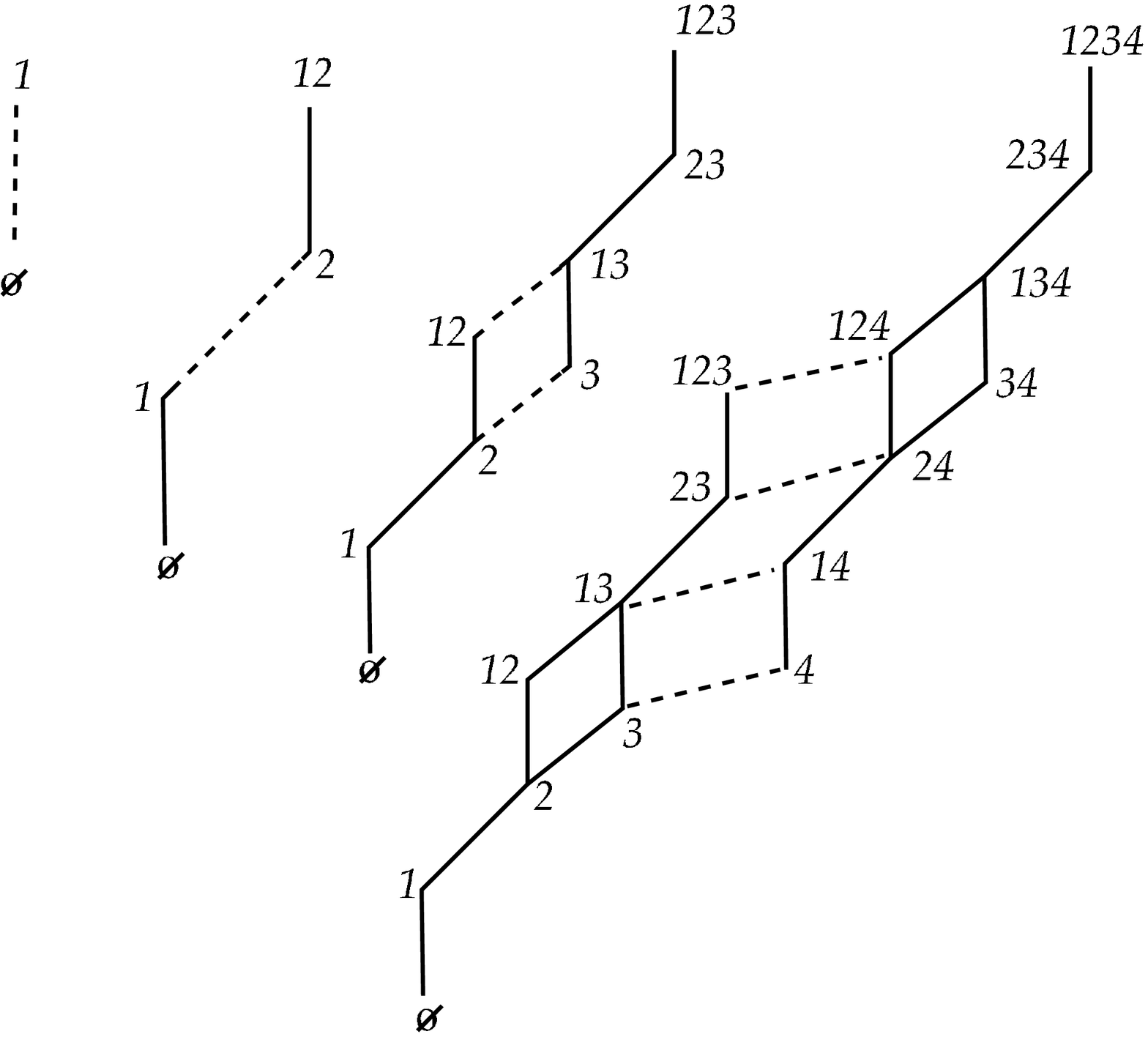, height=3.7in,width=4in}}
\caption{\sl The layered structure of $\eL_n$.}
\label{fig: mob1}
\end{figure}

\begin{ex} Let $n$ be a positive integer. Then  the poset $\eL_n$ has a natural layer structure $(\beps,\blam)$ defined by
\[
\beps(S) =\begin{cases}
1 & n\in S\\
-1 & n\not\in S.
\end{cases},\;\;\forall S\subset \{1,\dotsc, n\}, \;\;\mbox{and}\;\; \blam(S)= S\cup\{n\},\;\;\forall S\subset \{1,\dotsc, n-1\}.
\]
In Figure \ref{fig: mob1} we depicted the layer structures in the posets $\eL_n$, $1\leq n\leq 4$. The layers are separated  by dotted lines. \qed
\end{ex}

Suppose $P$ is a poset  equipped with a  layer structure $(\beps,\blam)$. We define the \emph{double} of $P$ to be the poset $\eD_P$ defined as  the connected sum
\[
\eD_P= P{}_{Q_0}\!\#_{Q_1} P
\]
where $Q= P^+$, $i_0: Q\ra P$ is the canonical inclusion $P^+\hra P$, and $i_1$ is the drop map $\bdel:P^+\ra P^-\hra P$. Equivalently,
\[
\eD_P=P\, {}_{P^+}\!\#_{P^-}\, P.
\]
As as set we can identify $\eD_P$ with the disjoint union $P\sqcup P=P\times \{-1,1\}$. We the identify $Q_0$ with $P^+\times \{-1\}$ and $Q_1$ with $P^-\times\{1\}$. The transition map $j_1: Q_0\ra Q_1$ is given by
\[
j_1(p,-1)= (\bdel(p), 1),\;\;\forall p\in P^+.
\]
We see that the double $\eD_P$ is equipped with a canonical layer structure $(\widehat{\beps}, \widehat{\blam})$, where
\[
\widehat{\beps}: P\times \{-1,1\}\ra \{-1,1\}
\]
 is the canonical projection and the lifting map $\widehat{\blam} : P\times \{-1\}\ra P\times \{1\}$ is the tautological map,
\[
\widehat{\blam}(x,-1)=(x,1),\;\;\forall x\in P.
\]
\begin{remark} Observe that if $(x_0,\eps_0),(x_1,\eps_1)\in \eD_P$ ($\eps_i=\pm 1$) then $(x_0,\eps_0)\prec_{\eD_P} (x_1,\eps_1)$ if and only if
\[
(\eps_0=\eps_1\;\;\mbox{and}\;\;x_0<_Px_1)\;\;\mbox{or}\;\;(\eps_0<\eps_1\;\;\mbox{and}\;\; \exists z\in P^+:\;\;x_0<_P z,\;\;\bdel(z)< x_1).\proofend
\]
\label{rem: comp}
\end{remark}
\begin{ex} Figure \ref{fig: mob1} shows that $\eL_{n+1}$ is the double of $\eL_n$ for $n=1,2,3$.\qed
\end{ex}

\begin{definition} Suppose $P$ is a poset   equipped with the layer structure $(\beps,\blam)$. We  say that the layer structure  satisfies  the property ($\bsM$) if for every $x\in P$ the sets $(P^+)^{\geq x}$ and $(P^-)^{\leq x}$ are nonempty  and  there exist maps
\[
\bm^-:P\ra P^- ,\;\;\bm^+:P\ra P^+
\]
such that
\[
\bm^+(x)=\min (P^+)^{\geq x},\;\;\bm^-(x)=\max (P^-)^{\leq x}.
\]
Equivalently,  this means that  if $x\in P^-$ and $y\in P^+$ then $x<y$ if and only if
\begin{equation}
x\leq \bm^-(y)\;\;\mbox{and}\;\;\bm^+(x)\leq y.
\label{eq: comm}
\end{equation}
We  say that $\bm^\pm$ are the \emph{associated maps} of the $\bsM$-structure.\qed
\end{definition}

\begin{lemma} Suppose $P$ is a poset equipped with a  layer structure $(\beps,\blam)$ satisfying property ($\bsM)$ with associated maps $\bm^\pm$. Then the double $\eD_P$ also satisfies the property ($\bsM$). Moreover, if we denote by $\widehat{\bm}^\pm:\eD_P\ra \eD_P^\pm$ the associated maps, and we identify $\eD_P$ with $P\times \{-1,1\}$  then
\[
\widehat{\bm}^+(x,\eps)= \begin{cases}
(x,1) &\beps=1,\\
(\,\bdel(\,\bm^+(x)\,), 1\,) &\eps=-1,
\end{cases},\;\;\whbm^-(x,\beps)=\begin{cases}
(x,-1) &\eps =-1,\\
(\,\blam(\,\bm^-(x)\,),-1\,) &\eps=1.
\end{cases}
\]
\label{lemma: mm}
\end{lemma}

\proof  Clearly $(\eD_P^+)^{\geq (x,1)}\neq\emptyset $ and in this case $\min (\eD_P^+)^{\geq (x,1)}=(x,1)$. Let  $(x,-1)\in \eD_P^-$.      Then
\[
x\leq_P\bm^+(x),\;\; (\bm_+(x),-1) \prec_{\eD} (\bdel(\bm^+(x)),1),
\]
so that
\[
(\bdel(\bm^+(x)),1)\in (\eD_P^+)^{\geq (x,-1)}.
\]
Suppose $(y,1)\in \eD_P^+$ and $(x,-1)\prec (y,1)$. Then (see Remark \ref{rem: comp}) there exists $z \in P^+$ and  such that
\[
x\leq_P z,\;\; \bdel(z) \leq_P  y.
\]
Hence $z\in (P^+)^{\geq x}$ so that $\bm^+(x)\leq z$, and thus $\bdel(\bm^+(x))\leq \bdel(z)\leq y$. Hence
\[
(\,\bdel(\bm^+(x)),\, 1\,)=\min (\eD_P^+)^{\geq (x,-1)}.
\]
Similarly, $(\eD_P^-)^{\leq (x,-1)}$ is nonempty and $\max (\eD_P^-)^{\leq (x,-1)}=(x,-1)$. Let $(x,1)\in \eD_P^+$. Then
\[
\bm^-(x)<_P x,\;\;(\blam(\bm^-(x)),-1\,)\prec_\eD  (x,1),
\]
so that
\[
(\blam(\bm^-(x)),-1\,)\in (\eD_P^-)^{\leq (x,1)}.
\]
Suppose $(y,-1)\leq  (x,1)$. Then, according to Remark \ref{rem: comp} there exists $z\in P^+$ such that
\[
y\leq_P z\;\;\mbox{and}\;\; \bdel(z) \leq_P x.
\]
Hence $\bdel(z) \in (P^-)^{\leq x}$ so that $\bdel(z)\leq \bm^-(x)$.   Using the monotonicity of $\blam$ we deduce
\[
y\leq_P z\leq_P\blam(\bm^-(x)),
\]
so that
\[
(\blam(\bm^-(x)),-1)=\max (\eD_P^-)^{\leq (x,1)}.\proofend
\]
Using Theorem \ref{th: mm} we deduce the following consequence.

\begin{corollary} Suppose $P$ is equipped with a layer  structure $(\beps,\blam)$ satisfying the property $\bsM$.  Denote by $\bmu$ the M\"{o}bius function of $P$, by $\bmu^+$ the M\"{o}bius function of $P^+$ and by $\widehat{\bmu}$ the M\"{o}bius function of the double $\eD_P$.  Then, for any $x_0,x_1\in P$, and any $\eps_0,\eps_1=\pm 1$ we have
\[
\widehat{\bmu}\bigl(\, (x_0,\eps_0),\;\;(x_1,\eps_1)\,\bigr)=\begin{cases}
\bmu(x_0,x_1) & \mbox{if}\;\;\eps_0=\eps_1,\\
-\bmu^+(x_0,\blam(x_1)) &\mbox{if}\;\;\eps_0<\eps_1,\;\;x_0\in P^+,\;\;x_1\in P^-,\\
0 & \mbox{otherwise}.
\end{cases}\proofend
\]
\label{cor: mmd}
\end{corollary}

\section{The  structure of $\eL_n$}
\label{s: 3}
\setcounter{equation}{0}
We want to apply the abstract  results proven so far to the
special case of the poset $\eL_n$.   We denote by $<_n$ the partial order on $\eL_n$, and by $\bmu_n$ the M\"{o}bius function of $\eL_n$. The following is the key structural result.

\begin{proposition}   For every $n\geq 1$   the map
\[
\Psi_n:\eD_{\eL_n}\ra \eL_{n+1},\;\;S\mapsto \Psi_n(S,\eps)=
\begin{cases}
S  &\eps=-1,\\
S\cup \{n+1\} & \eps=1
\end{cases}
\]
is an isomorphism of posets.
\label{prop: double}
\end{proposition}

\proof   We denote  by $\prec_n$ the partial order on $\eD_{L_n}$. Observe first that
\[
\eL_n^+:= \bigl\{\, S\in \eL_n;\;\;S\ni n\,\bigr\}.
\]
If $\bdel_n$ denotes the drop map of $\eL_n$    then
\[
\bdel_n(S)=S\setminus \{n\},\;\;\forall S\in \eL_n^+.
\]
Clearly $\Psi_n$ is a bijection.    We first prove that it is
increasing. Suppose
\[
(S_0,\eps_0)\prec_n(S_1,\eps_1).
\]
Using Remark \ref{rem: comp} we distinguish two cases.

\smallskip

\noindent $\bullet$ $\eps_0=\eps_1=\eps$.  Then $S_0<_{n}S_1$
which implies immediately that
$\Psi_n(S_0,\eps)<_{{n+1}}\Psi_n(S_1,\eps)$.

\noindent $\bullet$ $\eps_0=-1$, $\eps_1=1$ and  there exists
$T\in \eL_n^+$ such that
\[
S_0 <_{n} T,\;\; T\setminus\{n\} <_{n} S_1.
\]
Then
\[
S_0 <_n\bigl(\, T\setminus\{n\}\,\bigr)\cup\{n+1\} <_n S_1\cup
\{n+1\} \Longrightarrow \Psi_n(S_0,-1)<_{{n+1}}\Psi_n(S_1,1).
\]
This proves that $\Psi_n$ is increasing.    We have to prove that
the inverse map $\Phi_n=\Psi_n^{-1}$ is also increasing.   For a
set $S\in \eL_{n+1}$ we define $S'\in \eL_n$ by $S':=S\setminus \{n+1\}$. Observe  that
\[
\Phi_n^{-1}(S) = \begin{cases}
(S',-1) & (n+1)\not\in S\\
(S',1)& (n+1)\in S.
\end{cases}
\]
Suppose $S<_{n+1} T$. We distinguish  several cases.

\smallskip

\noindent $\bullet$  $(n+1)\in S$. Then since $S<_{n+1} T$ we must
also have $(n+1)\in T$ so that $S'<_n T'$ and thus
\[
(S',1)=\Phi_n(S)\prec_n (T',1)=\Phi_n(T).
\]
\noindent $\bullet$  $(n+1)\not\in S\cup T$.  In this case  again
we have  $S'=T'$ so that $S'<_n T'$ and thus
\[
(S',-1)=\Phi_n(S)\prec_n (T',-1)=\Phi_n(T).
\]
$\bullet$ $(n+1)\not\in S$, $(n+1)\in T$. Then
\[
S'=S,\;\;T'=T\setminus
\{n+1\},\;\;\Phi_n(S)=(S',-1),\;\;\Phi_n(T)= (T',1).
\]
To prove that  $\Phi_n(S)\prec_n \Phi_n(T)$  we need to find $U\in \eL_n^+$ such that
\[
S'\leq_n U,\;\;\bdel_n(U)=U\setminus \{n\} \leq _n T'.
\]
We define
\[
S'':= \begin{cases}
\emptyset & \mbox{if}\;\; S=\emptyset\\
S\setminus \{\max S \}& \mbox{if}\;\;S\neq\emptyset.
\end{cases}
\]
Note that  $S''< T'$.  Now  we set $U=S''\cup \{n\}\in\eL_n^+$.
Then
\[
S'\leq_n U\;\;\mbox{and}\;\;U\setminus \{n\} =S''<_n T'.\proofend
\]
From the above proposition we deduce inductively that we can
identify $\eL_n$ with the set of sequences $\vec{s}\in \{-1,1\}^n$
via the map
\[
\{-1,1\}^n\ni\vec{s}\mapsto X_{\vec{s}}=\bigl\{i;\;s_i=1\,\bigr\}.
\]
We then have
\[
\eL_n^\pm =\bigr\{ \vec{s};\;\;s_n=\pm 1\,\bigr\},
\]
The lifting map $\blam_n:\eL_n^-\ra \eL_n^+$ is given by
\[
(s_1,\dotsc, s_{n-1},-1)\mapsto (s_1,\dotsc, s_{n-1}, 1).
\]
We have natural isomorphisms $\vfi_\pm :\eL_{n-1}\ra \eL^\pm_n$
given by
\[
(s_1,\dotsc, s_{n-1})\mapsto (s_1,\dotsc, s_{n-1},\pm 1),
\]
and  we also have \emph{predecessor maps} $\bpi:\eL_n\ra \eL_{n-1}$ given by
\[
\bpi(s_1,\dotsc,s_{n-1},s_n)=(s_1,\dotsc,s_{n-1})\Longleftrightarrow \bpi(S)=S\setminus \{n\}.
\]
Note that  $\bpi\circ\vfi_\pm=\one$ and   $\vec{s}\leq_n\vec{t}$ if and only
\[
\sum_{i\geq k} s_i\leq \sum_{i\geq k} t_i,\;\;\forall k=1,\dotsc, n.
\]

 Using  Corollary \ref{cor: mmd}  and the above observations we deduce the following result.

\begin{corollary} For every $\vec{s},\vec{t}\in \eL_n$ we have
\[
\bmu_n\bigl(\,\vec{s}, \vec{t} \,\bigr)=\begin{cases}
\bmu_{n-1}\bigl(\,\bpi(\vec{s}),\bpi(\vec{t})\,\bigr)  & s_n=t_n\\
- \bmu_{n-2}\bigl(\,\bpi^2(\vec{s}),\bpi^2(\vec{t})\,\bigr) & s_n<t_n,\;\;s_{n-1}=1=-t_{n-1}.\\
0 & \mbox{otherwise}.
\end{cases}\proofend
\]
\label{cor: mob}
\end{corollary}

We can transform the above inductive formula  in  a  more explicit one. Given two sets $S,T\subset \{1,\dotsc, n\}$  we define
\[
\Delta_{S,T}:\{1,\dotsc, n\}\ra \bZ,\;\;\Delta_{S,T}(k) =\#(T\cap [k,n])-\#(S\cap [k,n]),
\]
The \emph{weight} of the pair $S,T$ is the integer
\[
w(S,T):= \sum_{k=1}^n \Delta_{S,T}(k).
\]
Note that
\[
S\leq_n T \Longleftrightarrow \Delta_{S,T}(k)\geq 0,\;\;\forall k=1,\dotsc,n.
\]
\begin{definition} We say that a pair $(S,T)$ of subsets of $\{1,\dotsc, n\}$ is \emph{elementary} if the following two conditions hold.

\begin{itemize}

\item $\Delta_{S,T}(k)\in \{0,1\}$, $\forall k=1,\dotsc, n$.

\item $ \Delta_{S,T}(k)\cdot \Delta_{S,T}(k+1)=0$, $\forall k=1,\dotsc, n-1$.
\end{itemize}\qed
\end{definition}

Thus, a pair $(S,T)\in \eL_n\times \eL_n$ is elementary if and only if in the sequence
\[
\Delta_{S,T}(1),\Delta_{S,T}(2),\dotsc, \Delta_{S,T}(n)
\]
we encounter only $0$'s and $1$'s,  but we do not encounter two consecutive  $1$'s. The weight of the elementary pair is then the number of $1$'s in the above sequence.

\begin{proposition} Let $S, T\subset  \{1,\dotsc, n\}$. Then the following statements are equivalent.

\smallskip

\noindent (a) The pair $(S,T)$ is elementary.

\noindent (b) There exists a  sequence
\[
1=\nu_1< \nu_2 <\cdots <\nu_k <n=\nu_{k+1},
\]
such that $\nu_j-\nu_{j-1}>1$, $\forall j=2,\dotsc, k$,  and (possibly empty) subsets
\[
A_1\subset B_1\subset\{1\}, \;\;C_j\subset (\nu_j+1, \nu_{j+1})\cap \bZ, \;\;j=1,\dotsc,k,
\]
 such that
\begin{subequations}
\begin{equation}
S= A_1\cup C_1\cup\{\nu_1\}\cup C_2\cup \cdots \{\nu_k\}\cup C_k,
\end{equation}
\begin{equation}
T=B_1\cup C_1\cup\{\nu_2+1\}\cup C_2\cup \cdots \cup \{\nu_k+1\}\cup C_k.
\label{eq: bead}
\end{equation}
\end{subequations}
\label{prop: elem}
\end{proposition}

\begin{remark} The technical  condition (b) can be easily visualized using the  beads-along-a-rod  picture      we described in the introduction.

We indicate a subset $T\subset \{1,\dotsc, n\}$ by placing  beads on a rod with   linearly ordered positions marked $1$ through $n$. An element $t\in T$ corresponds to a bed located in the position $t$ on the rod. Graphically, we depict by a ``$\bullet$'' the positions on the rod occupied by  a  bead, and by a ``$\circ$'' an  unoccupied position (see  Figure \ref{fig: mob2}).

We declare an element $t$ in $T$ to be mobile if either $t=1$ or $t-1\not\in T$. Graphically the mobile  elements  correspond to beads that can be moved one position to the left. (In the case of a bead located on the position $1$, moving it to the left corresponds to sliding it off the rod.) In Figure \ref{fig: mob2} the mobile positions are $1,4,7,10,12$.
\begin{figure}[h]
\centerline{\epsfig{figure=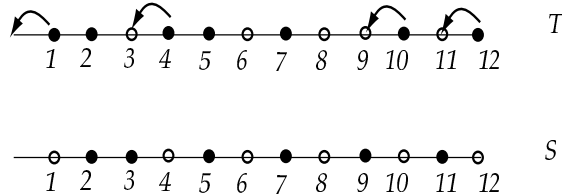, height=1in,width=2.9in}}
\caption{\sl Generating elementary pairs}
\label{fig: mob2}
\end{figure}

Condition (b)  signifies that the distribution  of beads $S$ is obtained from the distribution $T$ by sliding to the left  a certain number of mobile beads of $T$.  The number of beads that we slid to the left is precisely the weight of the pair $(S,T)$.       For the set $T$ described in (\ref{eq: bead}) the mobile  beads that are slid to the left are the beads located on the positions
\[
(B_1\setminus A_1)\cup\bigl\{ \,\nu_2+1,\dotsc,\nu_k+1\,\bigr\}.
\]
The weight of the pair $(S,T)$ is then $(k-1)+ \#(B_1\setminus A_1)$. In Figure \ref{fig: mob2} we obtain $S$ from $T$ by sliding to the left the mobile $T$-beads located at $1,4,10,12$. \qed
\label{rem: beads}
\end{remark}

The proof of Proposition \ref{prop: elem}  is an elementary induction  using the  above remark.  If we define
\[
\brho:\eL_n\ra \bZ,\;\;\brho(S):=\sum_{s\in S} s,\;\;\forall S\subset \{1,\dotsc,n\},
\]
then the sliding-beads description of elementary  pairs  implies the following result.

\begin{corollary}  If  $(S,T)$ is an elementary pair of subsets of $\{1,\dotsc, n\}$   then
\[
w(S,T)=\brho(T)-\brho(S).\proofend
\]
\label{cor: rank}
\end{corollary}

An elementary induction based on Corollary \ref{cor: mob} implies the following result.

\begin{theorem} Let $S, T\in \eL_n$. Then
\[
\bmu_n(S,T)=\begin{cases}
(-1)^{\brho(T)-\brho(S)} & \mbox{$(S,T)$ is  an elementary pair},\\
0 & \mbox{otherwise}.
\end{cases}
\tag{$\bmu$}
\]
\label{th: mob}
\end{theorem}

\begin{corollary}  Let $S,T\in \eL_n$.  If $\#T-\#S\not\in\{0,1\}$ then $\bmu_n(S,T)=0$.
\end{corollary}

\proof Observe that if $(S,T)$ is an elementary pair then $\#T-\#S\in \{0,1\}$.\qed

A chain of a poset $P$  is a linearly ordered  nonempty subset  $C\subset P$.    The endpoints of a chain $C$ are the elements $\min C$ and $\max C$.   The length of  a chain $C$ is   the integer
\[
\ell(C):=\# C -1.
\]
A poset  $P$ is called \emph{graded} if there exists an increasing function $r: P\ra  \bZ$  such that  for any  $x\leq y$ in $P$,  and any \emph{maximal} chain $C$ with endpoints $x$ and $y$  we have
\[
\ell(C) =r(y)-r(x).
\]
A function with this property is called a \emph{rank function} for the graded poset.

\begin{proposition} The poset $\eL_n$ is graded. As rank function we can take the function $\brho$.
\label{prop: graded}
\end{proposition}

\proof Clearly it suffices to prove that any  maximal  chain from $\emptyset$ to $S$ has length $\brho(S)$. Let $S\subset \{1,\dotsc, n\}$, $S\neq\emptyset$. A  maximal chain from $\emptyset$ to $S$  (of length $\ell$) is a sequence
\[
\emptyset \lessdot D_1\lessdot\cdots \lessdot  D_\ell= S
\]
where $I\lessdot J$      means that $J$ \emph{covers} $I$, i.e.,  $I<J$, and there is no element $K\in \eL_{n}$ such that $I<K<J$.  In terms  of  bead distribution  the condition $I\lessdot J$ signifies that the  bead distribution $I$ is obtained from $J$     after a single elementary left-slide. This shows that if $I\lessdot J$ then $\brho(J)=\brho(I)+1$ so that  every maximal chain from $\emptyset $ to $S$ has length $\rho(S)$. \qed

Recall that the dual of a poset $(P,<)$ is the poset $(P^*,<^*)$ which coincides with $P$ as a set but it is equipped with the opposite order, i.e.,
\[
x<^* y \Longleftrightarrow y<x.
\]
The M\"{o}bius function $\bmu^*$ of $P^*$ is related to the M\"{o}bius function of $P$ via the equality
\[
\bmu^*(x,y)=\bmu(y,x),\;\;\forall x,y\in P.
\]
A poset $P$ is called \emph{selfdual} if it is isomorphic to the dual poset $P^*$. Any poset isomorphism    $P\ra P^*$ is called a \emph{self-duality} of   $P$.

\begin{proposition} The map  $\bsi_n:\eL_n\ra \eL_n$ given by
\[
\{-1,1\}^n\ni \vec{s}\mapsto -\vec{s}\in \{-1,1\}^n
\]
is a selfduality of  $\eL_n$. In particular, we deduce that
\[
\bmu_n(\vec{s},\vec{t})=\bmu_n(-\vec{t},-\vec{s}).
\]
\end{proposition}

\proof  We have
\[
\vec{s}\leq\vec{t} \Longleftrightarrow \sum_{i\geq k} s_i <\sum_{i\geq k} t_i,\;\;\forall k=1,\dotsc, n
\]
\[
\Longleftrightarrow -\sum_{i\geq k} s_i >-\sum_{i\geq k} t_i,\;\;\forall k=1,\dotsc, n\Longleftrightarrow-\vec{t}<_n -\vec{t}.\proofend
\]
\begin{remark}   If  we regard an element $S\in \eL_n$ as a subset of $\{1,\dotsc, n\}$ then $\bsi_n(S)$ is   complement of $S$ in $\{1,\dotsc, n\}$.\qed
\label{rem: self}
\end{remark}

The selfduality $\bsi_n$  interacts nicely with the layer structure on $\eL_n$. More precisely, we have
\[
\bsi_n(\eL_n^\pm)=\eL_n^\mp.
\]
Lemma \ref{lemma: mm} and Proposition  \ref{prop: double} imply inductively that the layered poset $\eL_n$ satisfies the property ($\bsM$). We denote by $\bm_n^\pm:\eL_n\ra \eL_n^\pm$ the associated  functions.       Note that
\[
\bm_n^\pm\circ \bsi_n=\bsi_n\circ\bm_m^\mp.
\]
Using Lemma \ref{lemma: mm} we deduce the following result

\begin{lemma} Let $\vec{s}=(s_1,\dotsc,s_n)\in \{-1,1\}^n$, and set
\[
(t_1,\dotsc, t_{n-2},1)= \bm^+_{n-1}(s_1,\dotsc,s_{n-1}).
\]
 Then
\[
\bm_n^+(\vec{s})=\begin{cases}
\vec{s} & s_n=1\\
\bigl(\,\bdel_{n-1}\circ\bm_{n-1}^+(\bpi (\vec{s})\,), 1\,\bigr) & s_n=-1
\end{cases}
\]
\[
=\begin{cases}
\vec{s} & s_n=1\\
(t_1,\dotsc, t_{n-2}, -1, 1) & s_n=-1.
\end{cases}\proofend
\]
\label{lemma: mpl}
\end{lemma}

\begin{corollary} Let $S\subset \{1,\dotsc, n\}$ and set
\[
S':= \begin{cases}
\emptyset & \mbox{if}\;\; S=\emptyset\\
S\setminus \{\max S \}& \mbox{if}\;\;S\neq\emptyset.
\end{cases}
\]
Then
\[
\bm_n^+(S)= S'\cup\{n\}.
\]
\label{cor: mpl}
\end{corollary}

\proof   Define  $M_n: 2^{[n]}\ra 2^{[n]}$, $M_n(S)=S'\cup\{n\}$. It is easy to check that $\bm_1^+=M_1$ and that the maps  $M_n$ satisfy  the recurrence in Lemma \ref{lemma: mpl} so that $M_n=\bm_n^+$, $\forall n\geq 1$.\qed

\begin{remark} (a) The operation $S\mapsto \bm_n^+(S)$ has an intuitive description.  First $\bm_n^+(\emptyset)=\{n\}$.   If $n\in S$ then $\bm_n^+(S)=S$, while if $n\not\in S$, then $\bm_n^+(S)$ is obtained by ``trading'' the greatest element of $S$ for $n$, e.g.,
\[
\bm_{17}^+(\{ 2,5,7,8,\boldsymbol{11}\})=\{ 2,5,7,8,\boldsymbol{17}\}.
\]
In terms of bead distributions, the operation $S\mapsto \bm_n^+(S)$ corresponds to sliding the leftmost bead of $S$ all the way to the last position $n$  on the rod.

(b) Using the selfduality $\bsi_n$  we deduce
\[
\bm_n^-=\bsi_n\circ\bm_n^+\circ \bsi_n.
\]
Using    the description of $\bsi_n$ in Remark \ref{rem: self} we can give  very intuitive description of $\bm_n^-$. More precisely, if $n\not\in S$ then $\bm_n^-(S)=S$. Next, if $S=\{n\}$ then $\bm_n^-(S)=\emptyset$. Finally if $\{n\}\subsetneq S$ then $\bm_n^-(S)$ is obtained from $S$ by trading the element $n\in S$ with the greatest element \emph{not in $S$}. E.g.,
\[
\bm_{17}^-(\{ 2,5,7,8,16,\boldsymbol{17}\})=\{ 2,5,7,8,\boldsymbol{15},16\}.\proofend
\]
\label{rem: m-}
\end{remark}

We define  a \emph{simplicial scheme  with vertex set $V$}  to be a family  $\eS$ of nonempty subsets  $S\subset V$ such that
\[
\emptyset\neq S\subset T,\;\;T\in \eS\Longrightarrow S\in \eS.
\]
The sets $S\subset \eS$ are called the \emph{faces} of the simplicial scheme. To a simplicial scheme $\eS$ we can associate in a canonical way a triangulated space $|\eS|$ called the  \emph{geometric realization} of  $\eS$ (see \cite[\S 2]{GM}). Recall (see \cite{Bjorn}) that the \emph{nerve} of a poset $P$ is the simplicial scheme $\eN_P$  with vertex set $P$ defined by
\[
S\in \eN_P \Longleftrightarrow \mbox{$C$ is a chain in $P$}.
\]
We denote by $|P|$ the geometric realization of $\eN_P$.     We say that a poset $P$ is homeomorphic (homotopic) to a  topological space $X$ if its geometric realization is such.   The M\"{o}bius function $\bmu_P$ of $P$  is related to  $|P|$ via  the celebrated formula of P. Hall  (\cite[Eq. (9.14)]{Bjorn}, \cite[\S 3.8]{Sta})
\[
1+\bmu_P(x,y) = \chi\bigl(\, \bigl|\, (x,y)_P\,\bigr|\,\bigl),
\]
where $(x,y)_P$ denotes the open interval
\[
(x,y)_P:=\bigl\{ z\in P;\;\;  x<z<y\,\bigr\}
\]
and $\chi$ denotes the Euler characteristic of a space, with $\chi(\emptyset):= 0$.

\begin{theorem}  Suppose $(I,J)\in \eL_n\times \eL_n$ is an elementary pair. Then the closed interval  $[I,J]_{\eL_n}$ is isomorphic to the boolean poset $\eB_{\brho(J)-\brho(I)}$. In particular
\[
\bigl|\, (I,J)_{\eL_n}\,\bigr|\cong S^{\brho(J)-\brho(I)-2},
\]
where $S^k$ denotes the $k$-dimensional sphere if $k\geq 0$, while $S^{-1}:=\emptyset$.
\label{th: boole}
\end{theorem}

\proof  We argue by induction on $n$. The result is clearly true for $n=1$ so we assume it is true for any $k\leq n$ and we prove it for $n+1$.  Obviously, if $(n+1)\not\in J$ then $I, J\subset \{1,\dotsc, n\}$  and the  claim follows by induction.

Similarly, if $(n+1)\in I$ then $(n+1)\in J$,  the we have
\[
(I,J)_{\eL_{n+1}}\subset \eL_{n+1}^+\cong \eL_n,
\]
and again we can conclude by induction. Thus  we only need to consider the case $(n+1)\in J\setminus I$.  Since the pair $(I,J)$ is elementary we deduce from Proposition \ref{prop: elem} that $n\in S\setminus T$. We define
\[
\bar{I}:= I\setminus \{n\},\;\;\bar{J}:=  J\setminus \{n+1\}.
\]
Note that $\bar{I},\bar{J}\subset \{1,\dotsc, n-1\}$. From Proposition \ref{prop: elem} we deduce that $(\bar{I},\bar{J})\in \eL_{n-1}\times \eL_{n-1}$ is an elementary pair as well. Moreover
\begin{equation}
\brho(J)-\brho(I)=\brho(\bar{J})-\brho(\bar{I})+1.,
\label{eq: diff}
\end{equation}
and
\begin{equation}
\bm_{n+1}^+(I)=\bar{I}\cup\{n+1\},\;\;\bm_{n+1}^-(J)=\bar{J}\cup \{n\}
\label{eq: mpm}
\end{equation}
Observe that  if $S\in [I,J]_{\eL_{n+1}}$ then  $S$ \emph{contains exactly one of the  numbers $n$ or $n+1$}.    Now define a map
\[
[\bar{I},\bar{J}]_{\eL_{n-1}}\times \eB_1\cong [\bar{I},\bar{J}]_{\eL_{n-1}}\times \{ -1,1\} \stackrel{\Xi}{\Lra} [I,J]_{\eL_{n+1}},
\]
\[
[\bar{I},\bar{J}]_{\eL_{n-1}}\times \{ -1,1\}\ni (K, \eps)  \mapsto \Xi(K):= \begin{cases}
K\cup\{n\} &  \eps=-1\\
K\cup\{n+1\} &\eps=1.
\end{cases}
\]
We set
\[
\bar{S}:=S\setminus\{n,n+1\}\subset \{1,\dotsc, n-1\}.
\]
We distinguish two cases.

\smallskip

\noindent {\bf A.} \underline{$n\in S$}.  Then  $S=\bar{S}\cup \{n\}$. Using (\ref{eq: comm}) and the   inequality $S\leq_{n+1} J$  we deduce that
\[
S\leq \bm_{n+1}^-(J)=\bar{J}\cup \{n\}.
\]
Hence  $\bar{I}\cup \{n\} \leq _{n+1}S \leq_{n+1} \bar{J}\cup \{n\}$ so that $\bar{S}\in [\bar{I},\bar{J}]_{\eL_{n-1}}$.

\smallskip

\noindent {\bf B.} \underline{$n+1\in S$}.  Then $S=\bar{S}\cup\{n+1\}$. From the  inequality $I\leq_{n+1} S$ and (\ref{eq: comm}) we deduce that
\[
\bar{I}\cup \{n+1\} =\bm_{n+1}^+(I)\leq_{n+1} S=\bar{S}\cup \{n+1\}
\]
so that $\bar{S}\in [\bar{I},\bar{J}]_{\eL_{n-1}}$.     This discussion   shows that the map
\[
\Gamma: [I,J]_{\eL_{n+1}}\ra  [\bar{I}, {J}]_{\eL_{n-1}}
\]
given by
\[
[I,J]_{\eL_{n+1}}\ni S\mapsto \begin{cases}
(\bar{S},-1)  &  n\in S\\
(\bar{S}, 1) &n+1\in S.
\end{cases}
\]
is the inverse of the map  $\Xi$ so that $\Xi$ is a bijection. Clearly $\Xi$ is increasing.  Let us prove that $\Gamma$ is also increasing.  Suppose
\[
I\leq_{n+1} S<_{n+1} T\leq_{n+1} J.
\]
If $(n+1)\in S$ then $(n+1)\in T$  and we have $S=\bar{S}\cup \{n+1\} <_{n+1} T=\bar{T}=T\cup\{n+1\}$. Hence
\[
\Gamma(S)=(\bar{S}, 1)< (\bar{T},1)=\Gamma(T).
\]
We deduce similarly that if $n\in T$ then $\Gamma(S) <\Gamma(T)$. Thus we  need to discuss the case $n\in S$ and $(n+1)\in T$.    We have $S=\bar{S}\cup\{n\}$ and $T=\bar{T}\cup \{n+1\}$. Using  the inequality $S<_{n+1} T$ and (\ref{eq: comm}) we deduce
\[
\bar{S}\cup \{n\}\leq_{n+1}\bm_{n+1}^-(T)=\bar{T}\cup \{n\}\Longrightarrow \Gamma(S)=(\bar{S},-1)<(\bar{T},1)=\Gamma(T).
\]
Hence $\Xi$  is  an isomorphism of posets. From the induction assumption we deduce that we have an isomorphism of posets
\[
[I,J]_{\eL_{n+1}} \cong \eB_{\brho(\bar{J}) -\brho(\bar{I})}\times \eB_1,
\]
Using the isomorphism $\eB_{k+1}\cong \eB_k\times \eB_1$ and the equality  (\ref{eq: diff}) we  deduce that
\[
[I,J]_{\eL_{n+1}} \cong \eB_{\brho(J)-\brho(I)}.\proofend
\]

 \begin{proposition} The poset $\eL_n$ is a lattice. Moreover, for every $n>1$ and any  $\vec{s},\vec{t}\in \{-1,1\}^n$ such that $s_n\leq t_n$, we have
 \[
 \vec{s}\vee\vec{t} =\begin{cases}
 \bigl(\,\bpi(\vec{s})\vee \bpi(\vec{t}), \eps\,\bigr) & \mbox{  if $s_n=t_n=\eps\in\{\pm 1\}$}\\
\bigl(\,\bpi(\,\bm_n^+(\vec{s})\,)\vee \bpi(\vec{t}), 1\,\bigr) & \mbox{ if $s_n=-1$, $t_n=1$},
\end{cases}
\]
\[
 \vec{s}\wedge \vec{t} =\begin{cases}
 \bigl(\,\bpi(\vec{s})\wedge \bpi(\vec{t}), \eps\,\bigr) & \mbox{  if $s_n=t_n=\eps\in\{\pm 1\}$}\\
\bigl(\,\bpi(\vec{s}) \wedge \bpi(\,\bm_n^-(\vec{t})\,), -1\,\bigr) & \mbox{ if $s_n=-1$, $t_n=1$}.
\end{cases}
\]
\label{prop: lattice}
\end{proposition}

\proof  Using  the selfduality $\bsi_n$ it suffices  to prove  only   to deal with the join operation.  We argue by induction on $n$. Clearly $\eL_1$ is a lattice. For the inductive step consider $\vec{s},\vec{t}\in \eL_n$, $s_n\leq t_n$.  We distinguish two cases.

\smallskip

\noindent {\bf A.} If $s_n=t_n=\eps=\pm 1$ then $\vec{s},\vec{t}\in \eL_n^\pm$. Using the poset isomorphisms $\vfi_\pm:\eL_{n-1}\ra \eL_n^\pm$  and $\eL_n\cong \eD_{\eL_{n-1}}$ we deduce from the induction assumption that $\vec{s}\vee\vec{t}$ exists  and satisfies
\[
\vec{s}\vee\vec{t} =  \bigl(\,\bpi(\vec{s})\vee \bpi(\vec{t}), \eps\,\bigr).
\]
{\bf B.} $s_n=-1$, $t_n=1$. Then  $\vec{t}\in \eL_n^+$.    If $\vec{s},\vec{t}<_n\vec{v}$ then $\vec{v}\in \eL_n^+$. In particular, since $\eL_n$ satisfies property ($\bsM$) we deduce that
\[
\vec{v} {}_n\!>\bm_n^+(\vec{s}).
\]
Hence $\vec{s}\vee\vec{t}$ exists if and only if $\bm_n^+(\vec{s})\vee \vec{t}$ exists and in this case we have
\[
\vec{s}\vee\vec{t}=\bm_n^+(\vec{s})\vee \vec{t}.
\]
We  are now in the  case {\bf A} because  both $\bm_n^+(\vec{s})$ and  $\vec{t}$ belong to $\eL_n^+$.  The equality
\[
\bm_n^+(\vec{s})\vee \vec{t}= \bigl(\,\bpi(\,\bm_n^+(\vec{s})\,)\vee \bpi(\vec{t}), 1\,\bigr)
\]
follows as in case {\bf A}. \qed

\begin{ex}  Suppose $n>11$ and
\[
S=\{1,4,6,7,11\}\in \eL_n,\;\;T=\{2,5,9,10\}\in \eL_n.
\]
Then
\[
S\vee T= \{1,4,6,7,\boldsymbol{11}\}\vee \{2,5,9,\boldsymbol{11}\}= \{1,4,6,\boldsymbol{9},\boldsymbol{11}\}\vee \{2,5,\boldsymbol{9},\boldsymbol{11}\}
\]
\[
=\{1,4,\boldsymbol{6},\boldsymbol{9},\boldsymbol{11}\}\vee \{2,\boldsymbol{6},\boldsymbol{9},\boldsymbol{11}\}=\{1,\boldsymbol{4},\boldsymbol{6},\boldsymbol{9},\boldsymbol{11}\}\vee \{\boldsymbol{4},\boldsymbol{6},\boldsymbol{9},\boldsymbol{11}\}
\]
\[
=\{\boldsymbol{1},\boldsymbol{4},\boldsymbol{6},\boldsymbol{9},\boldsymbol{11}\}.
\]
\[
S\wedge T= \{1,4,6,7,\boldsymbol{10}\}\wedge \{2,5,9,\boldsymbol{10}\}=\{1,4,6,\boldsymbol{7},\boldsymbol{10}\}\wedge   \{2,5,\boldsymbol{7},\boldsymbol{10}\}
\]
\[
=\{1,4,\boldsymbol{5},\boldsymbol{7},\boldsymbol{10}\}\wedge   \{2,\boldsymbol{5},\boldsymbol{7},\boldsymbol{10}\}=\{1,\boldsymbol{2},\boldsymbol{5},\boldsymbol{7},\boldsymbol{10}\}\wedge   \{\boldsymbol{2},\boldsymbol{5},\boldsymbol{7},\boldsymbol{10}\}
\]
\[
= \{\boldsymbol{2},\boldsymbol{5},\boldsymbol{7},\boldsymbol{10}\}
\]
Observe that
\[
\brho(S\vee T)\brho(S\wedge T)=\brho(S)+\brho(T).\proofend
\]
\end{ex}

\begin{proposition} The  lattice  $\eL_n$  is modular, i.e.,
\begin{equation}
\brho(S\vee T)+\brho(S\wedge T)=\brho(S)+\brho(T),\;\;\forall S,T\in \eL_n.
\label{eq: modular}
\end{equation}
Moreover
\begin{equation}
S\vee\bsi_n(S)=\{1,\dotsc, n\},\;\;S\cap\bsi_n(S)=\emptyset,\;\;\forall S\in \eL_n,
\label{eq: comple}
\end{equation}
\begin{equation}
S\vee T=\{1,\dotsc, n\},\;\;S\wedge T=\emptyset \Longrightarrow T=\bsi_n(S).
\label{eq: comple1}
\end{equation}
\label{prop: modular}
\end{proposition}

\proof We argue by induction on $n$.  Denote by $\vee_n$ and $\wedge_n$ the lattice operations on $\eL_n$.

 For $n=1$ the result is obvious.   For the inductive step consider
\[
S=\{s_1<\dotsc <s_k\}\in\eL_{n+1},\;\;T=\{t_1<\cdots <t_\ell\}\in \eL_n,
\]
and  define $S'=S\setminus \{s_k\}$, $T':=T\setminus \{T_\ell\}$. Then Proposition \ref{prop: lattice} implies
\[
S\vee_{n+1} T =\bigl(\,S'\vee_n T'\,\bigr)\cup\{ \max(s_k,t_\ell)\}
\]
\[
S\wedge_{n+1} T =\Bigl(\,S'\wedge_n  T'\}\bigr)\cup\{ \min(s_k,t_\ell)\}.
\]
We deduce
\[
\rho(S\vee_{n+1} T)+\rho(S\wedge_{n+1}=\brho(S'\vee_n T')+ \brho(S\wedge_n T')+\max(s_k,t_\ell)+\min(s_k,t\ell)
\]
\[
\mbox{(by induction)}=\rho(S')+\brho(T') +s_k+t_\ell=\brho(S)+\brho(T).
\]
This proves  (\ref{eq: modular}). The inductive argument proving (\ref{eq: comple}) and (\ref{eq: comple1})  is very simple and we leave it to the reader. \qed

A modular lattice is $EL$-shellable . Using  \cite[Thm. 5.6]{BjWa}  (see also \cite[Sec. 3.2]{Wach}) we obtain the following result.

\begin{corollary} If $I<_n  J$ then  the geometric realization of the  open interval $(I,J)_{\eL_n}$ is contractible.\qed
\label{cor: contr}
\end{corollary}

We can be much more precise about the topology of the order intervals in the above corollary.

\begin{proposition} If $(I,J)$ is not an elementary pair of $\eL_n$ and $I<_n J$, then  the (open) order interval is homeomorphic to  $B^{\brho(J)-\brho(I)-2}$, where $B^d$ denotes the $d$-dimensional closed Euclidean ball.
\label{prop: ball}
\end{proposition}

\proof We argue as in the proof of  \cite[Thm. 2.7.7]{BB}. It suffices to investigate the structure of order intervals of length $2$.

\begin{figure}[h]
\centerline{\epsfig{figure=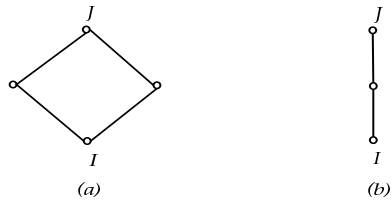, height=1in,width=2in}}
\caption{\sl Intervals of length $2$.}
\label{fig: mob4}
\end{figure}

Suppose that $I<_nJ$ and $\brho(J)-\brho(I)=2$. In terms of  distributions of beads along a rod this mens that the distribution $I$ can be obtained from the distribution $J$ by exactly two elementary left moves. If the two left moves involve different beads, so that the pair $(I,J)$ is elementary, then the Hasse diagram of $[I,J]_{\eL_n}$  is depicted in Figure \ref{fig: mob4}(a). If the two left moves involve the same bead, so that the pair $(I,J)$ is non-elementary, then $[I,J]_{\eL_n}$ is a chain of length 2 (Figure \ref{fig: mob4}(b)).

Using     Remark  \ref{rem: beads} we deduce that    if $I<_n J$ then the pair  $(I,J)$ is non elemntary if and only if there exists a non-elementary  pair $(I',J')$  such that  $I',J'\in [I,J]_{\eL_n}$, $I'<_nJ'$ and $\brho(J')-\brho(I')=2$.   This implies  (see \cite[Thm. 11.4]{Bjorn} or \cite[Sec. A2.4]{BB}) that the poset $(I,J)_{\eL_n}$ is homeomorphic to $B^{\brho(J)-\brho(I)-2}$.\qed

Recall that an element $x$ of a  lattice $P$ is called \emph{join reducible} if there exist $y,z<x$ such that $x=y\vee z$.

\begin{proposition} The set $S\subset \{1,\dotsc, n\}$ defines a join reducible element of the lattice $\eL_n$ if and only if  it has a \emph{gap}, i.e., there exists $1<k <n$ such that $k\not\in S$ and
\[
[1,k)\cap S, (k, n]\cap S\neq\emptyset.
\]
\label{prop: gap}
\end{proposition}

\proof  Suppose $S$ has a gap. Then $S$ has the form
\[
S=\{ s_1<\cdots <s_\ell\}\subset \{1,\dotsc, n\},
\]
 and for some $1< j\leq\ell$ we have $s_{j}-s_{j-1}>1$. We define
 \[
 S_0=\{ s_1,\dotsc, s_{j-1}, s_{j}-1, s_{j+1},\dotsc, s_\ell\},\;\;S_1= \{ s_{j}, s_j+1\dotsc,s_j+(\ell-j)\}.
 \]
 Then
 \[
 S_0, S_1<_n S,\;\;S_0\vee S_1=S.
 \]
 If  $S$  has no gap, so that $S$ has the form $S=\{j,j+1,\dotsc, j+\ell\}$, then $S$  covers a unique $S'\in \eL_n$, more precisely
 \[
 S'=\begin{cases}
 \{ j-1,j+1,\dotsc, j+\ell\} & j>1\\
 \{ j+1,\dotsc, j+\ell\} & j=1.
 \end{cases}
 \]
This proves that $S$   is  join irreducible because if $S_0,S_1<_n S$ then $S_0,S_1\leq_n S'$ so that $S_0\vee S_1\leq S'$. \qed


\end{document}